\documentclass[CRMATH,Unicode,manuscript]{cedram}

\usepackage{amsmath,amssymb}
\usepackage{amsthm}
\usepackage{latexsym}

\usepackage[english]{babel}

\usepackage[utf8]{inputenc}
\usepackage[T1]{fontenc}

\usepackage{color}
\usepackage{graphicx}

\newcommand{\dsp}{\displaystyle}

\newcommand{\al}{\alpha}

\newcommand{\ga}{\gamma}

\newcommand{\ld}{\lambda}

\newcommand{\si}{\sigma}

\newcommand{\eps}{\varepsilon}
\newcommand{\vphi}{\varphi}

\newcommand{\Rea}{\mathrm{Re}\,}

\def\R{\mathbb R}

\def\C{\mathbb C}

\renewcommand{\O}{\Omega}

\newcommand{\dx}{\mathrm{d}x}

\newcommand{\ds}{\mathrm{d}s}
\newcommand{\dt}{\mathrm{d}t}

\newcommand{\dv}{\mathrm{d}v}

\renewcommand{\do}{\mathrm{d}\omega}

\newcommand{\us}{\underset}

\newtheorem{theorem}{Theorem}
\newtheorem{corollary}{Corollary}
\newtheorem{lemma}{Lemma}
\newtheorem{definition}{Definition}
\newtheorem{proposition}{Proposition}

\newtheorem{counter-example}{Counter example}

\renewcommand\qedsymbol{$\blacksquare$}

\title{An elementary Tauberian proof of the Prime Number Theorem}

\author{\firstname{Philippe} \lastname{Angot}\CDRorcid{0000-0001-7535-0492}}
\address{Aix-Marseille Universit\'e, Institut de Math\'ematiques de Marseille,
CNRS UMR-7373 and Centrale Marseille, 3 place Victor Hugo, 13331 Cedex 03 Marseille, France}

\email[Ph. Angot]{philippe.angot@univ-amu.fr}



\begin{abstract}
We give a simple Tauberian proof\footnote{April 19, 2024} of the Prime Number Theorem using only elementary real analysis.
Hence, the analytic continuation of Riemann's zeta function $\zeta$ and its non-vanishing value on the whole line $\{z\in\C;\,\Rea z=1\}$ are no more required.
This is achieved by showing a strong extension for Laplace transforms on the real line of Wiener--Ikehara's theorem on Dirichlet's series, where the Tauberian assumption is reduced to a local boundary behavior around the pole.
\end{abstract}
\begin{document}

\maketitle

\section{Introduction and setting}
\label{sec:intro}

The prime\footnote{In the sequel, $p$ always denotes a prime number, {\em i.e.} a natural integer only divisible by $1$ and itself, and $p_n$ is the $n$-th prime number in the natural order.} number theorem (PNT), originally proved independently by J. Hadamard \cite{Hadamard_1896} and C. J. de la Vall\'ee Poussin \cite{Poussin_1896} in 1896 using ideas going back to B. Riemann (1859), see e.g. the monographs \cite{DKL_b2012,Korevaar_b2004}, reads as below:
\begin{equation}
\label{eq:PNT_pi_psi}
\dsp \pi(x) := \sum_{p\leq x} 1 \us{x\to+\infty}{\sim} \frac{x}{\ln x},
 \qquad\left(\mbox{ or equivalently: }\quad
 \psi(x) := \sum_{p^k\leq x} \ln p \us{x\to+\infty}{\sim} x\right).
\end{equation}
By introducing H. von Mangoldt's function (1895) $\Lambda$ as:
\begin{equation}
\label{eq:Mangoldt}
\dsp \Lambda(n) := \left\{
\begin{aligned}
\dsp &\ln p &\mbox{ if }\; n = p^{k} \quad\mbox{($p$ prime and integer $k\geq1$)}, \\
\dsp &0 &\mbox{ otherwise ($n$ is not a prime power)},
\end{aligned}
\right.
\end{equation}
the corresponding partial sum function is equal to the Chebyshev prime counting function $\psi$ defined by:
\begin{equation}
\label{eq:psi}
\dsp \psi(x) := \sum_{p^k\leq x} \ln p = \sum_{n\leq x} \Lambda(n),
 \qquad\forall x\geq1,
\end{equation}
summing over all prime powers not exceeding $x$.
It is relatively easy using the early estimates of P. L. Chebyshev (1850), see e.g. \cite{DKL_b2012}[Theorem 4.2], to show that the PNT (\ref{eq:PNT_pi_psi}) is actually equivalent to the claim:
\begin{equation}
\label{eq:PNT_psi}
\dsp \psi(x) \us{x\to+\infty}{\sim} x,
\end{equation}
that we are going to show.
It is also known from the Chebyshev estimates, e.g. \cite{DKL_b2012}[Theorem 2.6], that $\pi(x)=O(x/\ln x)$ or equivalently $\psi(x)=O(x)$ when $x\to+\infty$, and there is a constant $C>0$ such that:
\begin{equation}
\label{eq:psi_O}
\dsp \psi(x)
 = \sum_{n\leq x} \Lambda(n) \leq C\, x,
 \qquad\forall x\geq1.
\end{equation}

In this Note, we show that the present generalized Wiener--Ikehara's type Theorem \ref{theo:Ikehara_sdi} on the real line can be simply applied to prove in a quick way the Prime Number Theorem (\ref{eq:PNT_psi}).
Indeed, all the simple proofs of the PNT like Newman's Tauberian one (1980) \cite{Newman_1980} and its variants \cite{Korevaar_1982,Zagier_1997,Korevaar_2005}, except those of Erd\"os (1949) \cite{Erdos_1949}, Selberg (1949) \cite{Selberg_1949}, or their simplified version by Levinson (1969) \cite{Levinson_1969} or Daboussi (1984) \cite{Daboussi_1984}, use the crucial fact that the Riemann zeta function $\zeta$ has no zeros on the whole boundary line $\{z\in\C;\,\Rea z=1\}$.
It is thus an ultimate objective to seek a proof using as little about the zeta function as possible.
In the sequel, most of the technical requirements including the non-vanishing $\zeta$ on $\{\Rea z=1\}$ can be dropped since our generalization only requires local boundary behaviors on the real axis of the Laplace transform, and thus on $\zeta$. 
This makes the proof below far much easier and shorter.

\section{Tauberian theorems for Laplace transforms on the real axis}
\label{sec:Laplace_tauberian}
\subsection{Tauberian inversion of Abel--Laplace's summability}
\label{sec:inv_AL_summability}

The general setting in \cite{Angot_2023}, introduced for the Abel summability on power series, is extended in \cite{Angot_2024c} to the Abel--Laplace summability of complex-valued functions, where the Laplace transform is considered as a function of the real variable solely.
In particular, this enables us to deal with general or 'ordinary' Dirichlet's series.
Here, we restrict to real-valued functions.

\begin{definition}[\textbf{Slow decrease\footnote{The standard terminology is misleading since only the decrease of such functions is restricted, not their increase and so, any increasing function is slowly decreasing. This was introduced for sequences by R. Schmidt (1925); see \cite{Korevaar_b2004}.} at infinite}]
\label{def:sdi}
\hfill

A real-valued function $u:[0,+\infty[\rightarrow\R$ is said to be slowly decreasing (at infinite) if it satisfies the property:
\begin{equation}
\label{eq:sdi}
\dsp \lim_{\ld\to1^+} \liminf_{x\to+\infty} \big(u(\ld x)-u(x)\big) \geq 0.
\end{equation}
\end{definition}
\begin{proposition}[\textbf{Equivalent forms of the 'slow decrease' property}]
\label{prop:sdi}
\hfill

The 'slow decrease' property (\ref{eq:sdi}) of a function $u:[0,+\infty[\rightarrow\R$ is equivalent to the following one:
\begin{equation}
\label{eq:sdi_inf}
\dsp \lim_{\ld\to1^+} \liminf_{x\to+\infty}
 \left(\inf_{x\leq s\leq\ld x} \big(u(s)-u(x)\big) \right) \geq 0.
\end{equation}

Similarly, the 'slow decrease' property (\ref{eq:sdi}) of $u$ also reads equivalently as: 
\begin{equation}
\label{eq:sdi_sup}
\dsp \lim_{\ld\to1^+} \limsup_{x\to+\infty} \big(u(x)-u(\ld x)\big) \leq 0,
\end{equation}
that is equivalent to the following one:
\begin{equation}
\label{eq:sdi_limsup}
\dsp \lim_{\ld\to1^+} \limsup_{x\to+\infty}
 \left(\sup_{x\leq s\leq\ld x} \big(u(x)-u(s)\big) \right) \leq 0.
\end{equation}
\end{proposition}
\begin{small}
\begin{proof}
It is clear that (\ref{eq:sdi_inf}) implies (\ref{eq:sdi}). \\
Conversely, if (\ref{eq:sdi}) holds, then for all $\eps>0$, there exist $\eta>0$ and a sequence $(x_{n})$ tending to $+\infty$ such that for $1<\ld\leq 1+\eta$, we have:
\[
\dsp u(\ld x_{n}) - u(x_{n}) \geq -\eps,
 \qquad\mbox{ as }\; n\to+\infty.
\]
Hence, the infimum property ensures that:
\[
\dsp u(\ld x_{n}) - u(x_{n}) \geq \inf_{x_{n}\leq s\leq\ld x_{n}} \big(u(s)-u(x_{n})\big)
 \geq -\eps,
 \qquad\mbox{ as }\; n\to+\infty,
\]
which means that (\ref{eq:sdi_inf}) holds.

Moreover, (\ref{eq:sdi_sup}) is equivalent to (\ref{eq:sdi}) since $\limsup_{x\to+\infty}(-f(x))=-\liminf_{x\to+\infty}f(x)$ for any real-valued function $f$, and (\ref{eq:sdi_limsup}) is equivalent to (\ref{eq:sdi_sup}) similarly as above using the supremum property.
\end{proof}
\end{small}

We shall need the following preliminary results.
The first one is analogous to \cite{Angot_2023}[Lemma 2.4] on power series.

\begin{lemma}[\textbf{Tauberian result for Laplace transforms of positive functions}]
\label{lemma:tauberian_positive}
\hfill

Let $u:\R^{+}\to\R^+$ be a positive function in $L^{1}_{loc}([0,+\infty[)$ and consider the positive function $f:]0,+\infty[\to\R^{+}$ defined with its Laplace transform by:
\begin{equation}
\label{eq:F_Laplace}
\dsp f(t) := t \int_{0}^{+\infty} u(x) e^{-t x}\,\dx,
 \qquad\forall t>0,
\end{equation}
supposed to exist for all $t>0$ with a convergent integral.
We assume that:
\begin{equation}
\label{eq:f_A0}
\dsp f(0^+) := \lim_{t\to0^{+}} f(t) = 0.
\end{equation}
Then, we have:
\begin{equation}
\label{eq:C1_cv}
\dsp \lim_{x\to+\infty} \si(x)
 := \lim_{x\to+\infty} \frac{1}{x} \int_{0}^{x} u(s)\,\ds = 0.
\end{equation}
\end{lemma}
\begin{small}
\begin{proof}
We introduce the cut off function $k:[0,+\infty[\to\R^{+}$, positive and upper bounded, defined by:
\begin{equation}
\label{eq:k_cut-off}
\dsp k(x) := \left\{
\begin{aligned}
\dsp 0   &\quad\mbox{ for }\, 0\leq x < 1/e, \\
\dsp 1/x &\quad\mbox{ for }\, 1/e\leq x < +\infty, \\
\end{aligned}
\right.
\qquad\mbox{ such that: }\; 0\leq k(x)\leq e, \quad\forall x\in[0,+\infty[.
\end{equation}
Let us also consider the positive function $T:]0,+\infty[\to\R^{+}$ defined by:
\[
\dsp T(t) := t \int_{0}^{+\infty} u(x) e^{-t x} k(e^{-t x})\,\dx,
 \qquad\forall t>0.
\]
Then, with the positivity and the bound of $k$, we have:
\[
\dsp 0 \leq T(t) \leq e f(t),
 \qquad\forall t>0,
\]
and thus
\begin{equation}
\label{eq:T_bound}
\dsp 0 \leq \limsup_{t\to0^{+}} T(t) \leq e \limsup_{t\to0^{+}} f(t).
\end{equation}
Moreover, it results from the definition (\ref{eq:k_cut-off}) that for all $t>0$:
\[
\dsp e^{-t x} k(e^{-t x}) = 1, \quad\forall x\in[0,1/t],
 \qquad\mbox{ and }\quad
  k(e^{-t x}) = 0, \quad\forall x > 1/t.
\]
Thus, we get with such a truncature:
\begin{equation}
\label{eq:T_truncature}
\dsp T(t) = t \int_{0}^{1/t} u(x)\,\dx,
 \qquad\forall t>0.
\end{equation}
Now combining (\ref{eq:T_bound}) and (\ref{eq:T_truncature}) with the hypothesis (\ref{eq:f_A0}), it yields:
\begin{equation}
\label{eq:C1_cv_bis}
\dsp \lim_{t\to0^{+}} t \int_{0}^{1/t} u(s)\,\ds = 0,
\end{equation}
that is (\ref{eq:C1_cv}) by denoting $x=1/t$ so that $x\to+\infty$ when $t\to0^+$.
\end{proof}
\end{small}
\begin{lemma}[\textbf{Tauberian inversion of Ces\`aro summability for slowly decreasing functions}]
\label{lemma:Cesaro_inv_sdi}
\hfill

Let $u:\R^{+}\to\R$ be a slowly decreasing function in $L^{1}_{loc}([0,+\infty[)$ and define its Ces\`aro mean by:
\begin{equation}
\label{eq:Cesaro_mean}
\dsp \si(x) := \frac{1}{x} \int_{0}^{x} u(s) \ds
 = \int_{0}^{1} u(x t) \dt,
 \qquad\forall x>0, \qquad \si(0) := 0.
\end{equation}
We assume that:
\begin{equation}
\label{eq:C1_u}
\dsp \lim_{x\to+\infty} \si(x) = \ell\in\R.
\end{equation}
Then, we have:
\begin{equation}
\label{eq:C0_u}
\dsp \lim_{x\to+\infty} u(x) = \ell.
\end{equation}
\end{lemma}
\begin{small}
\begin{proof}
Let us observe with (\ref{eq:Cesaro_mean}) that:
\[
\dsp \si(x) - \ell = \frac{1}{x} \int_{0}^{x} (u(s)-\ell) \ds,
\]
and so, replacing $u$ by $u-\ell$ in (\ref{eq:Cesaro_mean}) if necessary, one may suppose that $\ell=0$ without loss of generality.

Now, let us notice that the assumption (\ref{eq:C1_u}) with $\ell=0$ implies that for any $\ld>1$: 
\begin{equation}
\label{eq:C1_u_mean}
\dsp \lim_{x\to+\infty} \frac{1}{\ld x-x} \int_{x}^{\ld x} u(s)\ds = 0,
 \qquad\forall \ld>1.
\end{equation}
Indeed, we have for all $x>0$ and any $\ld>1$:
\begin{equation*}
\begin{aligned}
\dsp \frac{1}{\ld x-x} \int_{x}^{\ld x} u(s)\ds
 &= \frac{1}{\ld x-x} \left(\int_{0}^{\ld x} u(s)\ds - \int_{0}^{x} u(s)\ds)\right) \\
 &= \frac{\ld}{\ld-1}\, \si(\ld x) - \frac{1}{\ld-1}\, \si(x)
  \quad\us{x\to+\infty}{\longrightarrow} 0, \qquad\forall \ld>1,
\end{aligned}
\end{equation*}
since, for any $\ld>1$, we have: $\si(\ld x)\to0$ when $x\to+\infty$, if $\si(x)\to0$ with (\ref{eq:C1_u}) and $\ell=0$.
Then, we have for any $\ld>1$ and a.e. $x>0$:
\begin{equation*}
\begin{aligned}
\dsp u(x)
 &= \frac{1}{\ld x-x} \int_{x}^{\ld x} u(s)\ds + \frac{1}{\ld x-x} \int_{x}^{\ld x} (u(x)-u(s))\ds \\
 &\leq \frac{1}{\ld x-x} \int_{x}^{\ld x} u(s)\ds + \sup_{x\leq s\leq\ld x} (u(x)-u(s)).
\end{aligned}
\end{equation*}
Now, passing to the upper limit when $x\to+\infty$ using (\ref{eq:C1_u_mean}), we get:
\[
\dsp \limsup_{x\to+\infty} u(x)
 \leq \limsup_{x\to+\infty} \left(\sup_{x\leq s\leq\ld x} (u(x)-u(s)) \right),
 \qquad\forall \ld>1,
\]
which gives by passing to the limit when $\ld\to1^+$ using the assumption of 'slow decrease' (\ref{eq:sdi_limsup}):
\begin{equation}
\label{eq:limsup_u}
\dsp \limsup_{x\to+\infty} u(x) \leq 0.
\end{equation}

Now, we proceed similarly as above by integrating in (\ref{eq:C1_u_mean}) over the interval $(x/\ld,x)$ instead of $(x,\ld x)$, and so we have for any $\ld>1$:
\begin{equation}
\label{eq:C1_u_mean2}
\dsp \lim_{x\to+\infty} \frac{\ld}{\ld x-x} \int_{x/\ld}^{x} u(s)\ds = 0,
 \qquad\forall \ld>1.
\end{equation}
Then, we have for any $\ld>1$ and a.e. $x>0$:
\begin{equation*}
\begin{aligned}
\dsp u(x)
 &= \frac{\ld}{\ld x-x} \int_{x/\ld}^{x} u(s)\ds
  + \frac{\ld}{\ld x-x} \int_{x/\ld}^{x} (u(x)-u(s))\ds \\
 &\geq \frac{\ld}{\ld x-x} \int_{x/\ld}^{x} u(s)\ds + \inf_{x/\ld\leq s\leq x} (u(x)-u(s)).
\end{aligned}
\end{equation*}
Now, passing to the lower limit when $x\to+\infty$ using (\ref{eq:C1_u_mean2}), we get:
\[
\dsp \liminf_{x\to+\infty} u(x)
 \geq \liminf_{x\to+\infty} \left(\inf_{x/\ld\leq s\leq x} (u(x)-u(s)) \right)
 = \liminf_{y\to+\infty}
 \left(\inf_{y\leq s\leq\ld y} \big(u(s)-u(y)\big) \right),
 \qquad\forall \ld>1,
\]
which gives by passing to the limit when $\ld\to1^+$ using the assumption of 'slow decrease' (\ref{eq:sdi_inf}):
\begin{equation}
\label{eq:liminf_u}
\dsp \liminf_{x\to+\infty} u(x) \geq 0.
\end{equation}
Hence, combining (\ref{eq:limsup_u}) and (\ref{eq:liminf_u}), it finally yields the desired result: $u(x)\to0$ when $x\to+\infty$.
\end{proof}
\end{small}

The following theorem can be improved by relaxing the one-sided Tauberian condition on $u$.
In particular, the lower bound $u\geq0$ is not necessary but this simplifies the proof.
Indeed, it is shown in \cite{Angot_2024c} that the 'slow oscillation' or 'slow decrease' Tauberian property, or also a still weaker one, is sufficient (and also necessary); see \cite{Angot_2023,Angot_2024} for similar results on power series. However, the present result below will be sufficient here for our purpose.

\begin{theorem}[\textbf{Tauberian inversion of the Abel--Laplace summability}]
\label{theo:AL_p-sdi}
\hfill

Let $u\in L^{1}_{loc}([0,+\infty[)$ be a positive and slowly decreasing real-valued function and
let us define the function $f:]0,+\infty[\to\C$ with the Laplace transform of $u$ by:
\begin{equation}
\label{eq:LT}
\dsp f(t) := t \int_{0}^{+\infty} u(x) e^{-t x}\,\dx
 = \lim_{X\to+\infty} t \int_{0}^{X} u(x) e^{-t x}\,\dx,
 \qquad\forall t>0,
\end{equation}
where the above improper integral is supposed to converge for any $t>0$.
We assume that for some finite value $\ell\in\R$:
\begin{equation}
\label{eq:AL0_f_u}
\dsp f(0^+) := \lim_{t\to0^{+}} f(t) = \ell\in\R.
\end{equation}
Then, we have:
\begin{equation}
\label{eq:C0_cv_u}
\dsp \lim_{x\to+\infty} u(x) = \ell.
\end{equation}
\end{theorem}
\begin{small}
\begin{proof}
The positive Abel--Laplace kernel: $k_{t}:x\mapsto e^{-t x}>0$, satisfies the unitary property:
\begin{equation}
\label{eq:AL_unit_kernel}
\dsp t \int_{0}^{+\infty} k_{t}(x)\dx = t \int_{0}^{+\infty} e^{-t x}\dx = 1,
 \qquad\forall t>0,
\end{equation}
and so, we have with (\ref{eq:LT}):
\[
\dsp f(t) - \ell = t \int_{0}^{+\infty} (u(x)-\ell) e^{-t x}\,\dx.
\]
Thus, whatever $u\in L^{1}_{loc}([0,+\infty[)$, by replacing $u(x)$ by $u(x)-\ell$ in (\ref{eq:LT}) if necessary and using (\ref{eq:AL_unit_kernel}), it suffices to study the case with $\ell=0$ without loss of generality.

Accordingly, let now $\ell=0$. With the positivity of $u$, Lemma \ref{lemma:tauberian_positive} ensures with (\ref{eq:AL0_f_u}) that $\si(x)\to0$ when $x\to+\infty$.
Then, with the 'slow decrease' (\ref{eq:sdi}) Tauberian assumption on $u$, Lemma \ref{lemma:Cesaro_inv_sdi} finally shows that: $u(x)\to0$ when $x\to+\infty$.
\end{proof}
\end{small}
\subsection{Extended Tauberian theorems of Wiener--Ikehara's type on the real axis}
\label{sec:Wiener-Ikehara}

The next Tauberian results on Dirichlet's series strongly extend the theorems of Wiener--Ikehara (1931)  \cite{Ikehara_1931} or Ingham's (1935) \cite{Ingham_1935}, originally proved using Wiener's (1932) Tauberian theory \cite{Wiener_1932} in Fourier analysis.
We refer to the review of complex Tauberian theory by Korevaar (2002) \cite{Korevaar_2002} or \cite{Korevaar_2005} and to the recent works of Debruyne--Vindas (2016 \& 2019) \cite{DV_2016,DV_2019} or Zhang (2019) \cite{Zhang_2019} for many further generalizations.
But all of them and their assumptions are concerning the complex plane.
On the contrary, the present generalization only requires local boundary assumptions on the real axis of the Laplace transform.
This considerably weakens the usual assumptions of analytic or continuous extension to the whole closed half-plane $\{z\in\C;\;\Rea z\geq1\}$ or also the $L^{1}_{loc}$-boundary behavior on the whole boundary line $\{z\in\C;\Rea z=1\}$ in the case of a simple pole at $z=1$; see \cite{Korevaar_2005,DV_2019} for details.
Hence, this strongly extends the previous results of this type.

\begin{theorem}[\textbf{Generalized Wiener--Ikehara's theorem on the real axis}]
\label{theo:Ikehara_sdi}
\hfill

Let $\al:[0,+\infty[\rightarrow\R^+$ be a positive function in $L^{1}_{loc}([0,+\infty[)$ and let us define its Laplace transform for some real number $\mu\geq0$:
\begin{equation}
\label{eq:LT_al_mu}
\dsp g(s) := \int_{0}^{+\infty} \al(x) e^{-s x}\dx
 := \lim_{X\to+\infty} \int_{0}^{X} \al(x) e^{-s x}\dx,
 \qquad\forall s>\mu\geq0,
\end{equation}
supposed to be well-defined for any $s>\mu\geq0$.
Next, we assume that there exist finite values $A,\ell\in\R$ such that:
\begin{equation}
\label{eq:L0_g_mu}
\dsp \lim_{s\to\mu^+} \left(g(s)-\frac{A}{s-\mu}\right) = \ell.
\end{equation}
Then (\ref{eq:L0_g_mu}) implies the following asymptotic behavior:
\begin{equation}
\label{eq:C0_al_mu}
\dsp \lim_{x\to+\infty} e^{-\mu x} \al(x) = A,
 \qquad\mbox{ i.e. }\quad \al(x)\us{x\to+\infty}{\sim} A\, e^{\mu x},
\end{equation}
{\em if} the positive function $\vphi:x\mapsto e^{-\mu x} \al(x)$ is slowly decreasing.
\end{theorem}
\begin{small}
\begin{proof}
The assumption (\ref{eq:L0_g_mu}) also reads:
\[
\dsp \lim_{s\to\mu^+} (s-\mu)\, g(s)
 = \lim_{s\to\mu^+} (s-\mu) \int_{0}^{+\infty} e^{-\mu x} \al(x) e^{-(s-\mu)x}\dx = A,
\]
which means defining the function $\vphi:x\mapsto e^{-\mu x} \al(x)$ on $[0,+\infty[$:
\begin{equation}
\label{eq:L0_phi_mu}
\dsp \lim_{t\to0^+} f(t) := \lim_{t\to0^+} t \int_{0}^{+\infty} \vphi(x) e^{-t x}\dx = A.
\end{equation}
Thus in these conditions, Theorem \ref{theo:AL_p-sdi} does apply for $\vphi$ and shows that:
\[
\dsp \lim_{x\to+\infty} \vphi(x) = A,
\]
which means that (\ref{eq:C0_al_mu}) holds true if the function $\vphi$ is slowly decreasing.
\end{proof}
\end{small}
\section{A simple proof of the Prime Number Theorem}
\label{sec:pnt}

Before entering to the details, let us show a useful property.

\begin{lemma}[\textbf{Slow decrease property with $\psi$}]
\label{lemma:psi_invx_sdi}
\hfill

The function $x\mapsto\psi(x)/x$ for $x>0$, where $\psi$ is the Chebyshev prime counting function in (\ref{eq:psi}), is slowly decreasing in the sense of Definition \ref{def:sdi}.
\end{lemma}
\begin{small}
\begin{proof}
For all $y,x>0$, we have:
\begin{equation}
\label{eq:pnt_psi_eq}
\begin{aligned}
\dsp \frac{\psi(y)}{y} - \frac{\psi(x)}{x}
 &= \frac{1}{y}\left(\psi(x) + (\psi(y)-\psi(x))\right) - \frac{\psi(x)}{x} \\
 &= - \frac{y-x}{y}\, \frac{\psi(x)}{x} + \frac{1}{y}\, (\psi(y)-\psi(x)).
\end{aligned}
\end{equation}
Since the function $\psi$ is monotonically increasing and using the Chebyshev bound (\ref{eq:psi_O}), it gives the lower bound below for all $y\geq x\geq1$:
\[
\dsp \frac{\psi(y)}{y} - \frac{\psi(x)}{x}
 \geq - C\, \left(1-\frac{x}{y}\right),
 \qquad\forall y\geq x\geq1.
\]
Now, passing to the lower limit when $x\to+\infty$ as $y/x\to1^+$, we obtain:
\begin{equation}
\label{eq:psi_invx_sdi}
\dsp \liminf_{x\to+\infty} \left(\frac{\psi(y)}{y} - \frac{\psi(x)}{x}\right) \geq 0,
 \qquad\mbox{ as }\quad 1<\frac{y}{x} \us{x\to+\infty}{\longrightarrow} 1^+,
\end{equation}
which means that the function $x\mapsto\psi(x)/x$ is slowly decreasing by satisfying (\ref{eq:sdi}).
\end{proof}
\end{small}

B. Riemann (1859) introduced the zeta function $\zeta$ with Dirichlet's series in the complex plane, extending the works of Euler and Dirichlet for the real variable. In the sequel, $\zeta$ is a function of the real variable $s>1$:
\begin{equation}
\label{eq:zeta}
\dsp \zeta(s) := \sum_{n=1}^{+\infty} \frac{1}{n^s}, \qquad\forall s>1.
\end{equation}
Now, all the properties that we need further on $\zeta:]1,+\infty[\rightarrow\R^+$ can be obtained by elementary real ana\-lysis, essentially to show that $(s-1)\,\zeta'(s)/\zeta(s)$ has a finite limit when $s\to1^+$.
In order to make the link with prime numbers, it is known from L. Euler (1737) by the unique representation of any positive integer $n$ as a product of prime powers, that $\zeta$ may be written as the product:
\begin{equation}
\label{eq:zeta_Euler}
\dsp \zeta(s) = \prod_{n=1}^{+\infty} \frac{1}{1-p_{n}^{-s}}
 = \left(\prod_{n=1}^{+\infty} \big(1-\frac{1}{p_{n}^{s}}\big)\right)^{-1},
 \qquad\forall s>1.
\end{equation}
Moreover, the comparison between series and integrals easily gives the bounds:
\begin{equation}
\label{eq:zeta_bound}
\dsp \frac{1}{s-1} \leq \zeta(s) \leq 1 + \frac{1}{s-1},
 \qquad\forall s>1.
\end{equation}
This means that $\zeta(s)$ is equivalent to $1/(s-1)$ when $s\to1^+$ and this simple pole at $s=1$ corresponds to the harmonic series.
Besides, since $\zeta(s)\geq1$, we get from (\ref{eq:zeta_bound}) that: $\zeta(s)\to1$ when $s\to+\infty$.
Still by series-integral comparison, the calculation for (\ref{eq:zeta_bound}) can be refined further to yield the first-order asymptotic expansion when $s\to1^+$:
\begin{equation}
\label{eq:zeta_asympt}
\dsp \zeta(s) = \frac{1}{s-1} + \ga + O(s-1), \quad\mbox{ as }\; s\to1^+,
\end{equation}
where $\ga$ is the Euler constant (1734).
Next, the logarithmic differentiation of the Euler product in (\ref{eq:zeta_Euler}) gives:
\begin{equation}
\label{eq:zeta_derivee1}
\dsp - \frac{\zeta'(s)}{\zeta(s)} = \sum_{n=1}^{+\infty} \frac{\Lambda(n)}{n^s},
 \qquad\forall s>1.
\end{equation}

Now, in our setting of Theorem \ref{theo:Ikehara_sdi}, we take:
\begin{equation}
\label{eq:al_psi}
\dsp \al(x) := \psi(e^{x}) = \sum_{n\leq e^{x}} \Lambda(n)
 = \sum_{\ln n\leq x} \Lambda(n),
 \qquad\forall x\geq0,
\end{equation}
in order to recover the integral representation of the Dirichlet's series given by (\ref{eq:zeta_derivee1}) where the partial sum function of coefficients $\psi$ is defined in (\ref{eq:psi}).
Indeed, it follows for $g$ defined in (\ref{eq:LT_al_mu}), first by the change of variable $v=e^{x}$, and then using Mellin transform of $\psi$ with $\psi(1)=0$:
\begin{equation}
\label{eq:g_al_psi}
\begin{aligned}
\dsp g(s) := \int_{0}^{+\infty} \al(x) e^{-s x}\dx
 &= \int_{0}^{+\infty} \psi(e^{x}) e^{-s x}\dx
  = \int_{1}^{+\infty} v^{-s-1}\psi(v) \dv \\
 &= - \frac{\zeta'(s)}{s\,\zeta(s)},
 \qquad\forall s>1.
\end{aligned}
\end{equation}
We refer to \cite{Korevaar_2002} or \cite{Korevaar_b2004}[Chap. III] for the details on the representation of a Dirichlet's series by the Mellin and Laplace transforms.
The equality (\ref{eq:g_al_psi}) can be also obtained directly, without introdu\-cing Mellin transform,  from the Dirichlet's series in (\ref{eq:zeta_derivee1}) with (\ref{eq:al_psi}) by using Abel's summation by parts with $n^{-s}=e^{-s \ln n}$, for all $n\geq1$.
Using now (\ref{eq:zeta_asympt}), it comes with (\ref{eq:g_al_psi}):
\begin{equation}
\label{eq:g_zeta_equiv}
\dsp g(s) \sim \frac{1}{s-1} \sim \zeta(s), \quad\mbox{ as }\; s\to1^+,
\end{equation}
and thus $g(s)$ behaves like $1/(s-1)$ as $\zeta(s)$ when $s\to1^+$.
More precisely, we deduce from (\ref{eq:zeta_asympt}) that:
\begin{equation}
\label{eq:g_zeta_asympt}
\dsp g(s) = \frac{1}{s(s-1)} - \frac{\ga}{s} + O(s-1)
 = \frac{1}{s-1} - \frac{1+\ga}{s} + O(s-1), \quad\mbox{ as }\; s\to1^+.
\end{equation}
Indeed, we get when $s\to1^+$:
\begin{equation*}
\begin{aligned}
\dsp \frac{\zeta'(s)}{\zeta(s)}
 &= \frac{\frac{-1}{s-1} + O(s-1)}{1+\ga(s-1)+O(|s-1|^2)}
  = \left(\frac{-1}{s-1} + O(s-1)\right) \left(1-\ga(s-1)+O(|s-1|^2)\right) \\
 &= \frac{-1}{s-1} + \ga + O(s-1), \quad\mbox{ as }\; s\to1^+.
\end{aligned}
\end{equation*}

Hence, it results from (\ref{eq:g_zeta_asympt}) that:
\begin{equation}
\label{eq:g_zeta_pole}
\dsp \lim_{s\to1^+} \left(g(s) - \frac{1}{s-1}\right) = - (1+\ga),
\end{equation}
and so the main assumption (\ref{eq:L0_g_mu}) in Theorem \ref{theo:Ikehara_sdi} is satisfied for $\ell=-(1+\ga)$ with $A=1$ and $\mu=1$.
Moreover, the positive function $\vphi:x\mapsto e^{-x}\al(x)=e^{-x}\psi(e^{x})$ is slowly decreasing as the function $x\mapsto\psi(x)/x$ from Lemma \ref{lemma:psi_invx_sdi}, since the exponential $x\mapsto e^{x}$ is increasing.
Hence, Theorem \ref{theo:Ikehara_sdi} ensures that: $\al(x)=\psi(e^{x})\sim e^{x}$ as $x\to+\infty$, and so $\psi(x)\sim x$ as $x\to+\infty$, which finally proves the PNT (\ref{eq:PNT_psi}).
\hfill $\blacksquare$

\medskip
By the way, we have the following equivalence result with PNT.

\begin{corollary}[\textbf{PNT-equivalent property on $\psi$}]
\label{cor:pnt_equiv_psi}
\hfill

The prime number theorem (\ref{eq:PNT_psi}) is equivalent to the slow decrease property of the function $x\mapsto\psi(x)/x$ in Lemma \ref{lemma:psi_invx_sdi}.

More precisely, the prime number theorem (\ref{eq:PNT_psi}) implies the following property on $\psi$:
\begin{equation}
\label{eq:pnt_psi}
\dsp \psi(y) - \psi(x) = (y-x) + o(y),
 \qquad\mbox{ as }\quad y\geq x \to+\infty.
\end{equation}
\end{corollary}
\begin{small}
\begin{proof}
It is clear using the PNT result (\ref{eq:PNT_psi}) that the function $x\mapsto\psi(x)/x$ satisfies (\ref{eq:sdi}) and thus, is slowly decreasing.

Conversely, if the function $x\mapsto\psi(x)/x$ is slowly decreasing as stated in Lemma \ref{lemma:psi_invx_sdi}, then the above proof shows the PNT result (\ref{eq:PNT_psi}).

Moreover, by passing to the limit when $y\geq x\to+\infty$ in Eq. (\ref{eq:pnt_psi_eq}) using (\ref{eq:PNT_psi}), we get:
\begin{equation}
\label{eq:pnt_psi_equiv}
\dsp \lim_{y\geq x\to+\infty} \frac{1}{y} \left(\psi(y) - \psi(x) - (y-x)\right) = 0,
\end{equation}
which means that (\ref{eq:pnt_psi}) holds true.
\end{proof}
\end{small}
\section{Conclusion}
\label{sec:conclusion}

As a consequence of the present result, because this is equivalent to PNT (\ref{eq:PNT_pi_psi}) from Hadamard (1896) \cite{Hadamard_1896} (see e.g. \cite{Korevaar_b2004}[Chap. III]) and since we have established here a proof of PNT independent of that, the Riemann zeta function $\zeta$ is completely free of zeros on the line $\{z\in\C;\,\Rea z=1\}$.

In conclusion, the new result is that it is possible to prove the PNT by a Tauberian method without assuming the non-vanishing of the Riemann function $\zeta$ on the whole line $\{z\in\C;\,\Rea z=1\}$ and  using only elementary real analysis, which was not clear at all until now.




\bigskip

\noindent\textbf{Acknowledgements}
\hfill

I would like to thank my colleague Olivier Ramar\'e (D.R. at CNRS, Institut de Math\'ematiques de Marseille) for the fruitfull discussions that have contributed to improve the first draft of this manuscript.

\bibliographystyle{crplain}
\bibliography{crmath_pnt_psd}

\end{document}